\documentclass[letterpaper, 10 pt, conference]{ieeeconf} 
\IEEEoverridecommandlockouts                              
\let\labelindent\relax
 \usepackage{etoolbox}
\makeatletter
\patchcmd{\@begintheorem}{\textit}{\textbf}{}{}
\makeatother
 \newtheorem{definition}{\bf Definition}

  \newtheorem{thm}{\bf Theorem}

 \newtheorem{lemma}{\bf Lemma}
  
  \usepackage{dirtytalk}
\usepackage{amsmath}
\usepackage{amsfonts}
\usepackage{amssymb}
\usepackage{bbm}
\usepackage{framed}
\usepackage{breqn}
\usepackage{subcaption}
\usepackage{cite}
\usepackage{enumitem}
\usepackage{cases}
\usepackage{pgfplots}
\usepackage{tabu}
\usepackage{url}
\usepackage{multicol}
\usepackage{color}
\usepackage[bottom]{footmisc}
\usepackage{environ}
\usepackage{tikz}
\usetikzlibrary{arrows,decorations.markings}
\usetikzlibrary{decorations.pathreplacing}

\usetikzlibrary{automata,arrows,positioning,calc}
\title{\LARGE\bf Unpredictable Planning Under Partial Observability}
\author{Michael Hibbard, Yagiz Savas, Bo Wu, Takashi Tanaka, Ufuk Topcu \thanks{This work was supported in part by the grants AFRL FA9550-19-1-0169, ONR N00014-19-1-2054, and DARPA D19AP00004 \newline \indent \quad All authors are with the Department of Aerospace Engineering
and Engineering Mechanics, and the Oden Institute for Computational
Engineering and Sciences, University of Texas, Austin, 201 E 24th
St, Austin, TX 78712. email: {\tt\small $\{$mwhibbard, yagiz.savas, bwu3, ttanaka, utopcu$\}$@utexas.edu}}}

\begin{document}
\maketitle
\begin{abstract}  We study the problem of synthesizing a controller that maximizes the entropy of a partially observable Markov decision process (POMDP) subject to a constraint on the expected total reward. Such a controller minimizes the predictability of a decision-maker's trajectories while guaranteeing the completion of a task expressed by a reward function. First, we prove that a decision-maker with perfect observations can randomize its paths at least as well as a decision-maker with partial observations. Then, focusing on finite-state controllers, we recast the entropy maximization problem as a so-called parameter synthesis problem for a parametric Markov chain (pMC). We show that the maximum entropy of a POMDP is lower bounded by the maximum entropy of this pMC. Finally, we present an algorithm, based on a nonlinear optimization problem, to synthesize an FSC that locally maximizes the entropy of a POMDP over FSCs with the same number of memory states. In numerical examples, we demonstrate the proposed algorithm on motion planning scenarios.
\end{abstract}
\section{Introduction}

A partially observable Markov decision process (POMDP) models sequential decision-making in stochastic environments with imperfect information and nondeterministic choices \cite{Kaelbling:1998:PAP:1643275.1643301, Kochenderfer}. A controller, i.e., a decision rule based on the imperfect information, resolves the nondeterminism and induces a stochastic process. In this paper, we are interested in synthesizing a controller that induces a stochastic process with maximum entropy among those whose realizations accumulate an expected total reward above a given threshold.

Entropy measures the unpredictability of outcomes in a random variable \cite{Cover}. Following \cite{Biondi,savas2018entropy}, we quantify the unpredictability of realizations in a stochastic process by defining the entropy of the process as the joint entropy of a sequence of random variables. Intuitively, our objective is then to synthesize a controller that induces a process whose realizations accumulate rewards in a way that maximizes the unpredictability to an outside observer.

Based on the previous sequence of actions and observations, a controller for a POMDP specifies a probability distribution over action selection. In the first part of the paper, we rigorously prove that an agent with perfect observations can randomize its trajectories at least as much as an agent with partial observations. We do so by proving that the maximum entropy of a POMDP is upper bounded by the entropy of its corresponding fully observable counterpart.

A finite-state controller (FSC) for a POMDP specifies a probability distribution over actions for each of its memory states according to the most recent information received from the environment \cite{Poupart,meuleau1999solving}. In this regard, FSCs represent a subset of controllers which may, in general, utilize the whole information history. Following the results of \cite{junges2017permissive}, we show that the entropy maximization problem over FSCs can be recast as a so-called parameter synthesis problem for a parametric Markov chain (pMC) \cite{Cubuktepe:10.1007/978-3-030-01090-4_10, Parametric}, under certain assumptions on the memory state transitions of the FSC. We can then efficiently solve for an entropy maximizing FSC for a fixed number of memory states. We first show that the maximum entropy of a pMC induced from a POMDP by FSCs with deterministic memory transitions is a lower bound on the maximum entropy of the POMDP. We also show that by using a specific memory transition function for FSCs, one can monotonically increase the maximum entropy of the stochastic process induced from a POMDP by increasing the number of memory states in the FSC. Finally, we present an algorithm, based on a nonlinear optimization problem introduced in \cite{Cubuktepe:10.1007/978-3-030-01090-4_10}, to synthesize parameters that maximize the entropy of a pMC subject to expected reward constraints.

An application of the proposed methods is the synthesis of a controller for an autonomous agent carrying out a mission in an adversarial environment. In particular, if the agent's sensor measurements are noisy and the mission is defined in terms of a reward function, the synthesized controller leaks the minimum information about the agent's trajectories to an outside observer while guaranteeing the accumulation of an expected total reward above a desired threshold.

\textbf{Related Work.} A recent study \cite{savas2018entropy} showed that an entropy-maximizing controller for an MDP could be synthesized efficiently by solving a convex optimization problem. In POMDPs, entropy has often been used for active sensing applications \cite{kreucher2005sensor,roy2005finding,eidenberger2010active}, where an agent seeks to select actions that maximize its information gain from the environment. These applications differ from our own as we seek to maximize the entropy of the trajectories an agent follows rather than maximizing its knowledge of the environment.

In the reinforcement learning literature, the entropy of a controller has been used as a regularization term in an agent's objective to balance the trade-off between exploration and exploitation \cite{haarnoja2018composable}. As discussed in \cite{haarnoja2017reinforcement}, using a controller with high entropy, an agent can learn a greater variety of admissible methods to complete a task, leading to a greater robustness when subsequently fine-tuned to specific scenarios. In imitation learning \cite{ziebart2008maximum}, a controller with high entropy similarly yields greater robustness when the provided demonstrations are imperfect. Unlike the aforementioned work, here we aim to synthesize a controller that maximizes the entropy of the induced stochastic process, rather than synthesizing a controller with high entropy. 

Synthesizing a controller for a POMDP subject to expected total reward constraints is, in general, undecidable in the infinite horizon case and PSPACE-complete in the finite horizon case \cite{undecidability}. For this reason, we focus on the space of FSCs \cite{amato2010optimizing,chatterjee2004trading}, which require only the most recent fragment of the actions and observations made by the agent. In doing so, we obtain locally optimal controllers over the space of FSCs with fixed numbers of memory states. We also provide a specific memory transition function that is guaranteed to increase the entropy of an induced stochastic process with respect to increasing number of memory states. 

\section{Preliminaries}
 For a set $\mathcal{S}$,  we denote its power set and cardinality by $2^{\mathcal{S}}$ and $\lvert \mathcal{S} \rvert$, respectively. The set of all probability distributions on a finite set $\mathcal{S}$, i.e., all functions $f $$:$$ \mathcal{S}$$ \rightarrow $$[0,1]$ such that $\sum_{s\in \mathcal{S}} f(s)$$=$$1$, is denoted by $\Delta(\mathcal{S})$. If $\{x_t\}$ is a sequence, a subsequence $(x_k,x_{k+1},\ldots,x_l)$ is denoted by $x_k^l$. We also write $x^l$$:=$$(x_1,x_2,\ldots,x_l)$. Finally, $\mathbb{N}$$=$$\{1,2,\ldots\}$, $\mathbb{N}_0$$=$$\{0,1,2,\ldots\}$ and $\mathbb{R}_{\geq 0}$$=$$[0,\infty)$.

\subsection{Partially Observable Markov Decision Processes}

{
\begin{definition}\label{def:MDP}
	A \textit{partially observable Markov decision process} (POMDP) is a tuple $\mathcal{M}=(\mathcal{S},s_I,\mathcal{A},\mathcal{P},\mathcal{Z},\mathcal{O},\mathcal{R})$ where $\mathcal{S}$ is a finite set of states, $s_I$$\in$$\mathcal{S}$ is a unique initial state, $\mathcal{A}$ is a finite set of actions, $\mathcal{P}$$:$$\mathcal{S}$$\times$$\mathcal{A}$$\rightarrow$$\Delta(\mathcal{S})$ is a transition function, $\mathcal{Z}$ is a finite set of observations, $\mathcal{O}$$:$$\mathcal{S}$$\rightarrow$$\Delta(\mathcal{Z})$ is an observation function, and $\mathcal{R}$$:$$\mathcal{S}$$\times$$\mathcal{A}$$\rightarrow$$\mathbb{R}_{\geq 0}$ is a reward function.
\end{definition}}
For simplicity, we assume that all actions $a$$\in$$\mathcal{A}$ are available in all states $s$$\in$$\mathcal{S}$. Additionally, we assume that only a single observation is available from the initial state, i.e., $\lvert \mathcal{O}(s_I)\rvert$$=$$1$. For the ease of notation, we denote the transition probability $\mathcal{P}(s'|s,a)$ and the observation probability $\mathcal{O}(z|s)$ by $\mathcal{P}_{s,a,s'}$ and $\mathcal{O}_{s,z}$, respectively. 

For a POMDP $\mathcal{M}$, we obtain the \textit{corresponding fully observable} MDP $\mathcal{M}_{fo}$ by setting $\mathcal{Z}$$=$$\mathcal{S}$ and $\mathcal{O}_{s,s}$$=$$1$ for all $s$$\in$$\mathcal{S}$. A Markov chain (MC) is a fully observable MDP where $\lvert\mathcal{A}\rvert $$=$$1$.

A \textit{system history} of length $t$$\in$$\mathbb{N}$ for a POMDP $\mathcal{M}$ is a sequence $h^t$$=$$(s_I,a_1,s_2,a_2,s_3,\ldots,s_t)$ of states and actions such that $\mathcal{P}_{s_k,a_k,s_{k+1}}$$>$$0$ for all $k$$\geq$$1$. We denote the set of all system histories of length $t$ by $\mathcal{H}^{t}$ and define the set of all system histories as $\mathcal{H}$$:=$$\cup_{t\in\mathbb{N}}\mathcal{H}^{t}$. Any system history $h^t$$=$$(s_I,a_1,s_2,\ldots,s_t)$ of length $t$ has an associated \textit{observation history} $o^t$$=$$(z_{I},a_1,z_{2},\ldots,z_{t})$ of length $t$$\in$$\mathbb{N}$. In general, there are multiple observation histories that are admissible for a given system history. We denote the collection of all observation histories of length $t$ by $Obs_{\mathcal{M}}^{t}$ and define the set of all observation histories as $Obs_{\mathcal{M}}$$:=$$\cup_{t\in \mathbb{N}}Obs_{\mathcal{M}}^{t}$. 

\noindent
\begin{definition}
A \textit{controller} $\pi$ for a POMDP $\mathcal{M}$ is a mapping $\pi$$:$$Obs_{\mathcal{M}}$$\rightarrow$$\Delta(\mathcal{A})$. Denote the set of all controllers by $\Pi(\mathcal{M})$.
\end{definition}

The probability that a controller $\pi$ takes the action $a$$\in$$\mathcal{A}$ for the observation history $o^t$$\in$$Obs_{\mathcal{M}}^{t}$ is denoted by $\pi(a| o^t)$. 

In general, a controller $\pi$$\in$$\Pi(\mathcal{M})$ may require the use of the entire observation history, which can be of an arbitrary length \cite{ross2014introduction}. By restricting controllers to use only the most recent fragment of their observation history, we obtain the class of controllers known as finite-state controllers \cite{ Meuleau, junges2017permissive}.

{
\begin{definition}\label{def:FSC}
 For a POMDP $\mathcal{M}$, a $k$-\textit{finite-state controller} ($k$-FSC) is a tuple $\mathcal{C}$$=$$(Q,q_1,\gamma,\delta)$, where $Q$$=$$\{q_1,q_2,\ldots,q_{k}\}$ is a finite set of memory states, $q_{1}$$\in$$Q$ is the initial memory state, $\gamma$$:$$Q$$\times$$\mathcal{Z}$$\rightarrow$$ \Delta(\mathcal{A})$ is a decision function and $\delta$$:$$Q$$\times$$\mathcal{Z}$$\times$$\mathcal{A}$$\rightarrow$$\Delta(Q)$ is a memory transition function. We denote the collection of all $k$-FSCs by $\mathcal{F}_k(\mathcal{M})$.
\end{definition}}

For a memory state $q$$\in$$Q$ of a k-FSC $\mathcal{C}$, we denote its set of successor memory states $q'$$\in$$Q$ by $Succ(q)$$:=$$\{q'$$\in$$ Q | \sum_{z\in\mathcal{Z}}\sum_{a\in\mathcal{A}}\delta(q'|q,z,a)$$>$$0\}$. 

 {
\begin{definition}\label{deterministicFSC}
    A \textit{deterministic $k$-FSC} $\mathcal{C}$$=$$(Q,q_1,\gamma,\delta)$ is a $k$-FSC such that for all $q$$\in$$Q$, $|Succ(q)| = 1$. We denote the collection of all deterministic k-FSCs by $\mathcal{F}_{k}^{det}(\mathcal{M})$.
\end{definition}}

An FSC prescribes a probability distribution for both the action selection $\gamma$ and the memory state update $\delta$ based on the most recent observation and the FSC's current memory state. 

\subsection{Entropy of Stochastic Processes}
The \textit{entropy of a random variable} $X$ with a countable support $\mathcal{X}$ and probability mass function (pmf) $p(x)$ is
\begin{align}
    H(X):=-\sum_{x\in\mathcal{X}}p(x)\log p(x).
\end{align}

We use the convention that $0$$\log$$0$$=$$0$. Let $(X_1,X_2)$ be a pair of random variables with the joint pmf $p(x_1,x_{2})$ and the support $\mathcal{X}\times \mathcal{X}$. The \textit{joint entropy} of $(X_1,X_2)$ is 
\begin{align}
\label{joint_entropy}
H(X_1,X_2):= -\sum_{x_1\in \mathcal{X}}\sum_{x_{2}\in \mathcal{X}}p(x_1,x_{2})\log p(x_1,x_{2}),
\end{align}\noindent
and the \textit{conditional entropy} of $X_2$ given $X_1$ is
\begin{align}
\label{conditional_entropy}
&H(X_2 | X_1):=-\sum_{x_1\in \mathcal{X}}\sum_{x_{2}\in \mathcal{X}}p(x_1,x_{2})\log p(x_2 |x_1).
\end{align}\noindent
The definitions of the joint and conditional entropies extend to collections of $k$ random variables as shown in \cite{Cover}. 
A discrete \textit{stochastic process} $\mathbb{X}$ is a discrete time-indexed sequence of random variables, i.e., $\mathbb{X}$$=$$\{X_k$$\in$$\mathcal{X}$ $:$ $k$$\in$$\mathbb{N}\}$. 
{
\noindent\begin{definition} (Entropy of a stochastic process) \cite{Biondi_thesis}
The \textit{entropy of a stochastic process} $\mathbb{X}$ is defined as 
 \begin{align}\label{entropy_def_stochastic}
 H(\mathbb{X}) :=\lim_{k\rightarrow \infty}H & ( X^{k}).
 \end{align}
 \end{definition}}

Recall that $X^k$$:=$$(X_1,X_2,\ldots,X_k)$. The above definition is different from the \textit{entropy rate} of a stochastic process, which is defined as $\lim_{k\rightarrow \infty}\frac{1}{k}H( X^k)$ when the limit exists \cite{Cover}. The limit in \eqref{entropy_def_stochastic} either converges to a non-negative real number or diverges to positive infinity \cite{Biondi_thesis}. 

For a POMDP $\mathcal{M}$, a controller $\pi$$\in$$\Pi(\mathcal{M})$ induces a discrete stochastic process $\{X_k$$\in$$\mathcal{S}$ $:$ $ k$$\in$$\mathbb{N}\}$ where each $X_k$ is a random variable over the state space $\mathcal{S}$. We denote the entropy of $\mathcal{M}$ under a controller $\pi$$\in$$\Pi(\mathcal{M})$ by $H^{\pi}(\mathcal{M})$. 
 {
\noindent \begin{definition}
(Maximum entropy of a POMDP) The \textit{maximum entropy of a POMDP} $\mathcal{M}$ is defined as
\begin{align}\label{max_ent_definition}
H(\mathcal{M}):=\sup_{\pi\in\Pi(\mathcal{M})}H^{\pi}(\mathcal{M}).
\end{align}
\end{definition}}

\section{Problem Statement}
We consider an \textit{agent} whose decision-making process is modeled as a POMDP and an \textit{outside observer} whose objective is to infer the states occupied by the agent in the future from the states occupied in the past. Being aware of the observer's objective, the agent aims to synthesize a controller that minimizes the predictability of its future states while ensuring that the expected total reward it collects exceeds a specified threshold. 

We measure the predictability of the agent's future states by the entropy of the underlying stochastic process. The rationale behind this choice can be better understood by recalling (see, e.g., \cite{Cover}) that for any given $n$$\in$$\mathbb{N}$ and $k$$\leq$$n$,
\begin{align}
    H(X^n)=&H(X^n_k|X^{k-1})+H(X^{k-1}).\label{derive_traj}
\end{align}
Therefore, by maximizing the value of the left hand side of \eqref{derive_traj}, one maximizes the entropy of the all future sequences $(X_k,\ldots,X_n)$ for any history of sequences $(X_{1},\ldots,X_{k-1})$. 

\noindent\textbf{Problem 1 (Constrained entropy maximization):} For a POMDP $\mathcal{M}$ and a constant $\Gamma$, synthesize a controller $\pi^{\star}$$\in$$\Pi(\mathcal{M})$ that solves
\begin{subequations}
\begin{align}\label{objective1}
    &\underset{\pi\in\Pi(\mathcal{M})}{\text{maximize}} \quad H^{\pi}(\mathcal{M})\\ \label{constraint1}
    &\text{subject to:} \quad \mathbb{E}^{\pi}\Big[\sum_{t=1}^{\infty}\mathcal{R}(S_t,A_t)\Big]\geq \Gamma.
\end{align}
\end{subequations}

By \cite{astrom}, for a reward function $\mathcal{R}: \mathcal{S} \times \mathcal{A}\rightarrow \mathbb{R}$,
\begin{align*}
    \sup_{\pi \in \Pi(\mathcal{M})}\mathbb{E}^{\pi}\Big[\sum_{t=1}^{N}\mathcal{R}(S_t,A_t)\Big]\leq \sup_{\pi \in \Pi(\mathcal{M}_{fo})}\mathbb{E}^{\pi}\Big[\sum_{t=1}^{N}\mathcal{R}(S_t,A_t)\Big].
\end{align*}
This inequality implies that an agent with perfect observations can collect an expected total reward at least as high as the expected total reward collected by an agent with imperfect observations. Since the objective function in the entropy maximization problem is quite different from the classical expected total reward objective, it is not obvious whether a similar claim holds for the entropy maximization problem. In the next section, we establish that an agent with perfect observations can indeed randomize its trajectories at least as well as an agent with imperfect observations.

It is known that deciding the existence of a controller satisfying constraint (\ref{constraint1}) is PSPACE-complete \cite{papadimitriou}. Therefore, the synthesis of globally-optimal controllers for Problem 1 is, in general, intractable. For this reason, in Section \ref{sec:reformUsingFSC}, we shift our focus to the subset of controllers known as \textit{finite state controllers} (FSCs). For a fixed number of memory states, we show that FSCs yielding locally optimal solutions to Problem 1 can be synthesized efficiently.

\section{An Upper Bound on the Maximum Entropy}

In this section, we establish that an agent with perfect observations can randomize its actions at least as well as an agent with imperfect observations. Formally, we show that
\begin{equation*}
    H^{\pi}(\mathcal{M}) \leq H^{\pi}(\mathcal{M}_{fo}).
\end{equation*}
Recall that for a POMDP $\mathcal{M}$, a controller $\pi$$\in$$\Pi(\mathcal{M})$ induces a stochastic process $\{X_k$$\in$$\mathcal{S}$ $:$ $ k$$\in$$\mathbb{N}\}$ whose entropy is
\begin{align}\label{eq:defEntropyMDP}
   H^{\pi}(\mathcal{M}):=\lim_{k \rightarrow \infty} H^{\pi}(X^k) =\sum_{t=2}^{\infty} H^{\pi}(X_{t}|X^{t-1}),
\end{align}
where $H(X_1)$$=$$0$ because $\mathcal{M}$ has a unique initial state.

For a given system history $h^t$$=$$(s_I,a_1,s_2,a_2,s_3,\ldots,s_t)$, let the sequences $s^t$$=$$(s_I,s_2,s_3,\ldots,s_t)$ and $a^t$$=$$(a_1,a_2,a_3,\ldots,a_t)$ be the corresponding state and action histories of length $t$, respectively. We denote the set of all state and action histories of length $t$ by $\mathcal{SH}^t$ and $\mathcal{AH}^t$. Additionally, we define the set of all possible state and action histories as $\mathcal{SH}$$:=$$\cup_{t\in \mathbb{N}} \mathcal{SH}^t$ and $\mathcal{AH}$$:=$$\cup_{t\in \mathbb{N}} \mathcal{AH}^t$. 

For a POMDP $\mathcal{M}$ under the controller $\pi$$\in$$\Pi(\mathcal{M})$, it can be shown that the realization probability $Pr^{\pi}(s^{t+1} | s^{t})$ of the state history $s^{t+1}$$\in$$\mathcal{SH}^{t+1}$ for a given $s^t$$\in$$\mathcal{SH}^{t}$ is
\begin{align}
    Pr^{\pi}(s^{t+1} | s^{t})=\sum_{a^t\in\mathcal{AH}^t} \prod_{k=1}^t \mu_k(a_k| h^k) \mathcal{P}_{s_t,a_t,s_{t+1}}
\end{align}
where $h^k$ are prefixes of $h^t$ from which the state sequence $s^t$ is obtained, and $\mu_t$ $:$ $\mathcal{H}^t$$\rightarrow$$\Delta(\mathcal{A})$ is a mapping such that
\begin{align}\label{MDP_controller}
   \mu_t(a|h^t):= \sum_{o^{t} \in Obs_{\mathcal{M}}^{t}} \pi(a|o^{t}) Pr(o^{t}|h^{t})
\end{align}
where the realization probability $Pr(o^{t}|h^{t})$ of the observation history $o^{t}$ for a given $h^{t}$ can be recursively written as
\begin{align}\label{eq:obsSeqRecursive}
    Pr(o^{t}|h^{t}) = \mathcal{O}_{s_{t},z_{t}}\mathcal{P}_{s_{t-1},a_{t-1},s_{t}}Pr(o^{t-1}|h^{t-1})
\end{align}
\noindent
for all $t$$>$$1$ by assuming that $o^1$$=$$s_I$ with probability 1. 

Now, for a given controller $\pi$$\in$$\Pi(\mathcal{M})$ and a finite constant $T$$\in$$\mathbb{N}$, let $\mathcal{V}_{t,T}^{\pi}$ $:$ $\mathcal{SH}^t$$\rightarrow$$\mathbb{R}$ be the \textit{value function} such that
\begin{align}\label{eq:defEntropyMDPsum}
    \mathcal{V}_{t,T}^{\pi}(s^{t}) := \sum_{k=t}^{T}H^{\pi}(X_{k+1}|X^k_{t},X^{t}=s^{t}).
\end{align}

{
\noindent \begin{lemma}\label{writeValueFunction} For a POMDP $\mathcal{M}$, a controller $\pi$$\in$$\Pi(\mathcal{M})$ and a finite constant $T$$\in$$\mathbb{N}$,
    \begin{align}\label{eq:EntropyVfRecursive}\hspace{-0.3cm}
        \mathcal{V}_{t,T}^{\pi}(s^{t}) =& H^{\pi}(X_{t+1}|X^{t}=s^{t}) \\ \nonumber
        & + \sum_{s^{t+1} \in \mathcal{SH}^{t+1}} Pr^{\pi}(s^{t+1}| s^{t}) \mathcal{V}_{t+1,T}^{\pi}(s^{t+1}) \nonumber
    \end{align}
for all $t$$<$$T$ and $s^{t} \in \mathcal{SH}^{t}$.
\end{lemma}}
\textbf{Proof:} See Appendix.$\Box$

\noindent It is worth noting that 
\begin{align}\label{entropy_equivalence}
  \sup_{\pi\in\Pi(\mathcal{M})} H^{\pi}(\mathcal{M})=\sup_{\pi\in\Pi(\mathcal{M})}\lim_{T\rightarrow \infty} \mathcal{V}_{1,T}^{\pi}(s_I).
\end{align}
Moreover, since $\mathcal{V}^{\pi}_{t,T}$ is monotonically increasing in $T$ for all $\pi$$\in$$\Pi(\mathcal{M})$, we have, for all $s^t$$\in$$\mathcal{SH}^t$,
\begin{align}\label{limsup}
    \sup_{\pi\in\Pi(\mathcal{M})}\lim_{T\rightarrow \infty} \mathcal{V}_{t,T}^{\pi}(s^t)=\lim_{T\rightarrow \infty} \sup_{\pi\in\Pi(\mathcal{M})}\mathcal{V}_{t,T}^{\pi}(s^t).
\end{align}

As a consequence of Lemma \ref{writeValueFunction}, we can now define functions $\mathcal{V}^{\star}_{t,T}$ $:$ $\mathcal{SH}^t$$\rightarrow$$\mathbb{R}$  for $t$$\leq$$T$ such that
\begin{align}
    \mathcal{V}_{t,T}^{\star}(s^{t}) := \sup_{\pi\in\Pi(\mathcal{M})}\mathcal{V}_{t,T}^{\pi}(s^{t})
\end{align}
and conclude that, for all $t$$<$$T$ and $s^{t} \in \mathcal{SH}^{t}$, 
    \begin{align}\hspace{-0.3cm}\label{bellman_result}
        \mathcal{V}_{t,T}^{\star}(s^{t}) =&\sup_{\pi\in\Pi(\mathcal{M})} \Big[H^{\pi}(X_{t+1}|X^{t}=s^{t}) \\ \nonumber
        & + \sum_{s^{t+1} \in \mathcal{SH}^{t+1}} Pr^{\pi}(s^{t+1}| s^{t}) \mathcal{V}_{t+1,T}^{\star}(s^{t+1})\Big]. \nonumber
    \end{align}
Using \eqref{entropy_equivalence}, \eqref{limsup}, and taking the limit of both sides of (\ref{bellman_result}) as $T$$\rightarrow$$\infty$, we conclude that $H(\mathcal{M})$$=$$\lim_{T\rightarrow \infty}\mathcal{V}^{\star}_{1,T}(s_I)$ satisfies the equations in \eqref{bellman_result} which are recursive Bellman equations \cite{puterman2014markov}. 

Recall that for any controller $\pi$$\in$$\Pi(\mathcal{M})$ on a POMDP $\mathcal{M}$, we can construct, using \eqref{MDP_controller}, a controller $\pi'$$\in$$\Pi(\mathcal{M}_{fo})$ on the corresponding MDP $\mathcal{M}_{fo}$ which satisfies $Pr^{\pi}(s^{t+1} | s^{t})$$=$$Pr^{\pi'}(s^{t+1} | s^{t})$ for all $s^t,s^{t+1}$$\in$$\mathcal{SH}$. Then, 
\begin{align*}
    \sup_{\pi\in\Pi(\mathcal{M})}H^{\pi}(X_{t+1}|X^{t}=s^{t})\leq \sup_{\pi\in\Pi(\mathcal{M}_{fo})}H^{\pi}(X_{t+1}|X^{t}=s^{t})
\end{align*}

\noindent
for all $s^t$$\in$$\mathcal{SH}$. Informally, by having access to the state history $s^t$, a controller $\pi'$$\in$$\Pi(\mathcal{M}_{fo})$ can achieve an immediate reward $H^{\pi'}(X_{t+1}|X^{t}$$=$$s^{t})$ in \eqref{bellman_result} that is at least as high as the immediate reward achieved by a controller $\pi$$\in$$\Pi(\mathcal{M})$. We can then conclude the following result.
{
\noindent\begin{thm}\label{thm:POMDPbounded}
    For a POMDP $\mathcal{M}$ and its corresponding fully observable MDP $\mathcal{M}_{fo}$, we have 
\end{thm}}
\begin{equation}\label{eq:POMDPbounded}
        H(\mathcal{M}) \leq H(\mathcal{M}_{fo}).
\end{equation}
\textbf{Proof:} See Appendix.$\Box$

Based on the result of Theorem 1, we see that an agent with perfect observations can randomize its trajectories at least as well as an agent with partial observations.

\section{Reformulation using Finite-State Controllers}\label{sec:reformUsingFSC}

Since the synthesis problem over general controllers is, in general, intractable due to constraint (\ref{constraint1}), in this section, we consider the entropy maximization problem over deterministic FSCs with fixed numbers of memory states.

\noindent\textbf{Problem 2 (Constrained entropy maximization over FSCs):} For a POMDP $\mathcal{M}$ and constants $k$$>$$0$ and $\Gamma$, synthesize (if it exists) a controller $\mathcal{C}^{\star}$$\in$$\mathcal{F}^{det}_k(\mathcal{M})$ that solves
\begin{subequations}
\begin{align}\label{objective3}
    &\underset{\mathcal{C}\in\mathcal{F}^{det}_k(\mathcal{M})}{\text{maximize}}\ \ H^{\mathcal{C}}(\mathcal{M})\\ \label{constraint3}
    &\text{subject to:} \ \ \mathbb{E}^{\mathcal{C}}\Big[\sum_{t=1}^{\infty}\mathcal{R}(S_t,A_t)\Big]\geq \Gamma.
\end{align}
\end{subequations}

\subsection{A Solution Approach Through Parametric Markov Chains}
We develop solution methods to Problem 2 through the use of parametric Markov chains.
Recall that for a POMDP $\mathcal{M}$, a k-FSC $\mathcal{C}$$\in$$\mathcal{F}_k(\mathcal{M})$ induces a Markov chain (MC). The collection of all MCs that can be induced from $\mathcal{M}$ by a k-FSC is described by the induced parametric MC which is defined as follows. 
 {\noindent
\begin{definition}\label{def:pMC}
    For a POMDP $\mathcal{M}$ and a constant $k$$>$$0$, the \textit{induced parametric Markov chain} (pMC) is a tuple $\mathcal{D}_{\mathcal{M},k}$$=$$ (S_{\mathcal{M},k},s_{I,\mathcal{M},k},V_{\mathcal{M},k},P_{\mathcal{M},k})$ where
    \begin{itemize}
        \item $S_{\mathcal{M},k} = \mathcal{S} \times \{1,2,...,k\}$ is the finite set of states,
        \item $s_{I,\mathcal{M},k} = \langle s_{I}, 1 \rangle$ is the initial state,
        \item $\!\begin{aligned}[t]
                    V_{\mathcal{M},k} = \{\gamma_{a}^{q,z}| & z \in \mathcal{Z},q \in Q, a \in \mathcal{A} \} \nonumber \\
                    & \cup \{\delta_{q'}^{q,z,a}|z \in \mathcal{Z}, q,q' \in Q, a \in \mathcal{A} \} \nonumber
                \end{aligned}$\\
                is the finite set of parameters,
        \item $P_{\mathcal{M},k}$ $:$ $S_{\mathcal{M},k}$$\rightarrow$$\Delta(S_{\mathcal{M},k}) $ 
        is a transition function such that $P_{\mathcal{M},k}(s'|s):= \sum_{a\in A}\overline{P}(s'|s,a)$ for all $s,s'$$\in$$S_{\mathcal{M},k}$
 
    where $\overline{P}$ $:$ $S_{\mathcal{M},k}$$\times$$\mathcal{A}$$\rightarrow$$ \Delta(S_{\mathcal{M},k})$ is a mapping such that
     \end{itemize}
         \begin{equation}\label{parametrictransition}
        \overline{P}(\langle s',q' \rangle \ | \ \langle s,q \rangle,a) := \sum_{z \in \mathcal{Z}} \mathcal{O}_{s,z} \, \mathcal{P}_{s,a,s'} \, \gamma_{a}^{q,z} \, \delta_{q'}^{q,z,a}.
    \end{equation}
\end{definition}}
Note that when defining the (parametric) transition probabilities of the induced pMC, we suppose that the observations $\mathcal{O}_{s,z}$ are obtained before selecting actions $a$$\in$$\mathcal{A}$. We also remark that different definitions of the induced pMC can be used to reduce the number of parameters in $V_{\mathcal{M},k}$ \cite{junges2017permissive}. 

Now, an MC can be obtained from the induced pMC by instantiating the parameters $\mathcal{V}_{\mathcal{M},k}$ in a way that the resulting transition function $P_{\mathcal{M},k}$ is well-defined. Formally, let $Z$$=$$\{p_1,\ldots,p_n\}$ be a finite set of parameters over the domain $\mathbb{R}$, and $\mathbb{Q}[Z]$ be the set of multivariate polynomials over $Z$. An \textit{instantiation} for $Z$ is a function $u$$:$$ Z $$\rightarrow$$\mathbb{R}$. Additionally, replacing each parameter $p_i$ in a polynomial $f$$\in$$\mathbb{Q}[V]$ by $u(p_i)$ yields $f[u]$$\in$$\mathbb{R}$. 

Applying an instantiation $u$$:$$ V_{\mathcal{M},k} \rightarrow \mathbb{R}$ to the induced pMC $\mathcal{D}_{\mathcal{M},k}$, denoted $\mathcal{D}_{\mathcal{M},k}[u]$, replaces each polynomial $P_{\mathcal{M},k}$ by $P_{\mathcal{M},k}[u]$. An instantiation $u$ is then \textit{well-defined} for $\mathcal{D}_{\mathcal{M},k}$ if the replacement yields probability distributions, i.e., if $\mathcal{D}_{\mathcal{M},k}[u]$ is an MC. 

Every well-defined instantiation $u$ describes a k-FSC $\mathcal{C}_u$$\in$$\mathcal{F}_k(\mathcal{M})$ \cite{junges2017permissive}. Thus, we can synthesize all admissible MCs that can be induced from a POMDP $\mathcal{M}$ by a k-FSC $\mathcal{C}_{u}$$\in$$\mathcal{F}_k(\mathcal{M})$ through well-defined instantiations $u$ over $V_{\mathcal{M},k}$. This implies Problem 2 can be reduced to a parameter synthesis problem for the induced pMC. In Section \ref{FSCsynthesis}, for a pMC, we present a method to synthesize parameters that induces a stochastic process with maximum entropy whose realizations satisfy an expected total reward constraint. 

In the next section, we provide two results that allow one to compare the maximum entropy of a POMDP with the maximum entropy of the induced pMC. 

\subsection{An Upper Bound and a Monotonocity Result}
For a given $k$-FSC $\mathcal{C}$, let $u_{\mathcal{C}}$$:$$V_{\mathcal{M},k}$$\rightarrow$$\mathbb{R}$ be the corresponding instantiation of $\mathcal{D}_{\mathcal{M},k}$ such that $u_{\mathcal{C}}(\gamma_{a}^{q,z})$$:=$$\gamma(a|q,z)$ and $u_{\mathcal{C}}(\delta_{q'}^{q,z,a})$$:=$$\delta(q'|q,z,a)$. Note that $\mathcal{D}_{\mathcal{M},k}[u_{\mathcal{C}}]$ is a stochastic process. For a given POMDP $\mathcal{M}$ and a memory bound $k$$>$$0$, let the maximum entropy of the induced pMC $\mathcal{D}_{\mathcal{M},k}$ be defined as
\begin{align} H(\mathcal{D}_{\mathcal{M},k}):= \sup_{\mathcal{C} \in \mathcal{F}_{k}^{det}(\mathcal{M})}H(\mathcal{D}_{\mathcal{M},k}[u_{\mathcal{C}}]).
\end{align}
{\noindent\begin{thm}\label{thm:pMCbounded} Let $\mathcal{M}$ be a POMDP, $k$$>$$0$ be constant, and $\mathcal{D}_{\mathcal{M},k}$ be the induced pMC. Then,
\end{thm}}
\begin{equation}\label{eq:pMCbounded}
    H(\mathcal{D}_{\mathcal{M},k}) \leq H(\mathcal{M}).
\end{equation}
\textbf{Proof:} See Appendix.$\Box$

Theorem \ref{thm:pMCbounded} implies that by synthesizing a deterministic $k$-FSC $\mathcal{C}$ such that the instantiation $u_{\mathcal{C}}$ maximizes the entropy of the induced pMC $\mathcal{D}_{\mathcal{M},k}$, we can guarantee that the entropy $H^{\mathcal{C}}(\mathcal{M})$ of the POMDP $\mathcal{M}$ under the controller $\mathcal{C}$ is at least as high as the entropy of $H(\mathcal{D}_{\mathcal{M},k}[u_{\mathcal{C}}])$.

We now present a subclass of deterministic $k$-FSCs, for which we can monotonically increase the maximum entropy of a stochastic process induced from a POMDP by increasing the number of memory states in the FSC.  

For a POMDP $\mathcal{M}$, consider a $k$-FSC $\mathcal{C}$$=$$(Q,q_1,\gamma,\delta)$ with the memory transition function $\delta$$:$$Q$$\times$$\mathcal{Z}$$\times$$\mathcal{A}$$\rightarrow$$\Delta(Q)$
\begin{align}\label{LastLoopDetFSC}
\begin{cases}
        \delta(q_{i+1}|q_{i},z,a)=1 & \forall z \in \mathcal{Z}, a \in \mathcal{A}, 1 \leq i < k \\
        \delta(q_{k}|q_{k},z,a)=1  &\forall z \in \mathcal{Z}, a \in \mathcal{A}\\
        \delta(q_i|q_j,z,a)=0 & \text{otherwise}.
        \end{cases}
\end{align}

\begin{figure}[t]\centering
\scalebox{0.65}{
\begin{tikzpicture}[->, >=stealth', auto, semithick, node distance=2cm]

     \tikzstyle{every state}=[fill=white,draw=black,thick,text=black,scale=1]

     \node[state,initial,initial text=] (q_1) {\Large{$q_1$}};
     \node        (z_1) [below right=1mm and 3mm of q_1] {\Large{$z_1$}};
     \node        (z_2) [above right=1mm and 3mm of q_1] {\Large{$z_2$}};
     \node[state] (q_2) [right=20mm of q_1]  {\Large{$q_2$}};
     \node        (q_dots) [right=10mm of q_2] {\Large{$\cdots$}};
     \node[state] (q_n) [right=10mm of q_dots]  {\Large{$q_k$}};
     \node        (z_11) [below right=-1mm and 4mm of q_n] {\Large{$z_1$}};
     \node        (z_21) [above right=-1mm and 4mm of q_n] {\Large{$z_2$}};
    
    \draw[thick,-] (q_1) -- (z_1);
    \draw[thick,-] (q_1) -- (z_2);
    \draw[thick,-] (q_n) -- (z_11);
    \draw[thick,-] (q_n) -- (z_21);
    \path

    (z_1)  edge[out=320, in=230]   node[right]{\Large{$a_2$}}     (q_2)
    (z_1)  edge  node[below]{\Large{$a_1$}}     (q_2)
    (z_2)  edge   node[above]{\Large{$a_1$}}     (q_2)
    (z_2)  edge [out=40, in=130]   node[right]{\Large{$a_2$}}     (q_2)
    (q_2) edge (q_dots)
    (q_dots)  edge       (q_n)
       
    (z_11)  edge[out=225, in=250]   node[below]{\Large{$a_1,a_2$}}     (q_n)
    (z_21)  edge [out=135, in=110]   node[above]{\Large{$a_1,a_2$}}     (q_n);

\end{tikzpicture}}
\caption{ A deterministic $k$-FSC example. }\label{deterministic_controller_example}
\end{figure}
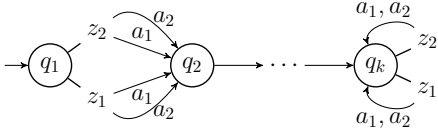

A $k$-FSC with the memory transition function defined above is shown in Fig. \ref{deterministic_controller_example}. Let $\overline{\mathcal{F}}_{k}(\mathcal{M})$$\subset$$\mathcal{F}_{k}^{det}(\mathcal{M})$ be the set of $k$-FSCs whose memory transition function is given in \eqref{LastLoopDetFSC}. Then, we have the following result.

{\begin{lemma}\label{lemma:pMCmonotonic}
The following inequality holds for all $j$$\leq$$k$.  
\end{lemma}}
\begin{equation}
    \sup_{\mathcal{C} \in \overline{\mathcal{F}}_{j}(\mathcal{M})}H(\mathcal{D}_{\mathcal{M},j}[u_{\mathcal{C}}])\leq \sup_{\mathcal{C} \in \overline{\mathcal{F}}_{k}(\mathcal{M})}H(\mathcal{D}_{\mathcal{M},k}[u_{\mathcal{C}}]).
\end{equation}
\textbf{Proof:} See Appendix.$\Box$

Based on the result of Lemma \ref{lemma:pMCmonotonic}, we can now set an initial number of memory states for a \textit{deterministic} FSC with the memory transition function \eqref{LastLoopDetFSC} and solve Problem 2 to determine the maximum entropy of the induced pMC. We may then iteratively adjust the number of memory states in the FSC to achieve a greater maximum entropy. 

\section{Finite-State Controller Synthesis}\label{FSCsynthesis}
 
We now present a method to synthesize a deterministic $k$-FSC that maximizes the entropy of a POMDP over all deterministic $k$-FSCs whose memory transition function is given in \eqref{LastLoopDetFSC}. 

Recall that for a POMDP $\mathcal{M}$ and a constant $k$$>$$0$, the induced pMC represents all possible MCs that can be induced from $\mathcal{M}$ by a $k$-FSC. Additionally, the maximum entropy of the induced pMC provides a lower bound on the maximum entropy of the POMDP due to Theorem \ref{thm:pMCbounded}. Furthermore, by increasing the number of memory states in $k$-FSCs with transition function given in \eqref{LastLoopDetFSC}, we can synthesize controllers that improves the entropy of the induced stochastic process.

Using Lemma \ref{writeValueFunction}, for a POMDP $\mathcal{M}$ and a constant $k$$>$$0$, we can write the entropy of an instantiation $u$$:$$V_{\mathcal{M},k}$$\rightarrow$$\mathbb{R}$ of the induced pMC $\mathcal{D}_{\mathcal{M},k}$, denoted $\mathcal{D}_{\mathcal{M},k}[u]$, as a solution to a form of Bellman equations. Specifically, let $T$$\subseteq$$S_{\mathcal{M},k}$ be the set of absorbing states in $\mathcal{D}_{\mathcal{M},k}[u]$, i.e., $s$$\in$$T$ implies that the only successor state of $s$ is itself. Let $P_{\mathcal{M},k}^{u}$$:$$S_{\mathcal{M},k}$$\rightarrow$$\Delta(S_{\mathcal{M},k})$ be the transition function of the instantiated pMC such that $P_{\mathcal{M},k}^{u}(s'|s)$ is defined by replacing parameters $\gamma_a^{q,z}$ and $\delta^{q,z,a}_{q'}$ in \eqref{parametrictransition} with their corresponding instantiations $u(\gamma_a^{q,z})$ and $u(\delta^{q,z,a}_{q'})$. Additionally, let $L^{u}$$:$$S_{\mathcal{M},k}$$\rightarrow$$\mathbb{R}$ be the \textit{local entropy} function such that
\begin{align}
L^{u}(s):=-\sum_{s'\in S_{\mathcal{M},k}}P_{\mathcal{M},k}^{u}(s'|s)\log P_{\mathcal{M},k}^{u}(s'|s)
\end{align}
\noindent
for all $s$$\in$$S_{\mathcal{M},k}$. Using Lemma \ref{writeValueFunction} and defining variables $\nu$$\in$$\mathbb{R}^{\lvert S_{\mathcal{M},k}\rvert}$, it can be shown that the entropy of $\mathcal{D}_{\mathcal{M},k}[u]$ is the unique fixed-point of the system of equations
\begin{subequations}
\begin{flalign}
  &  \nu(s) = L^{u}(s) + \sum_{s' \in S_{\mathcal{M},k}}P_{\mathcal{M},k}^{u} (s'| s)\nu(s')  \ s\in S_{\mathcal{M},k}\backslash T && \raisetag{23pt}\\
   & \nu(s)=0  \qquad \qquad \qquad \qquad \qquad \qquad \quad\ \  s\in T, && \label{fixed2}
    \end{flalign}
    \end{subequations}
such that $H(\mathcal{D}_{\mathcal{M},k}[u])$$=$$\nu(s_{I,\mathcal{M},k})$. Then, the maximum entropy $H(\mathcal{D}_{\mathcal{M},k})$ of $\mathcal{D}_{\mathcal{M},k}$ can be computed by finding the maximum $\nu(s_{I,\mathcal{M},k})$ that satisfies 
\begin{flalign}
  &  \nu(s) \leq L^{u}(s) + \sum_{s' \in S_{\mathcal{M},k}}P_{\mathcal{M},k}^{u} (s'| s)\nu(s')  \ s\in S_{\mathcal{M},k}\backslash T && \raisetag{23pt}
    \end{flalign}
together with the condition \eqref{fixed2}. Similarly, for the expected total reward constraint, let $\mathcal{R}^u$$:$$\mathcal{S}_{\mathcal{M},k}$$\rightarrow$$\mathbb{R}$ define the expected immediate rewards on $\mathcal{D}_{\mathcal{M},k}$ such that, for all $s$$\in$$S_{\mathcal{M},k}$,
\begin{align}
    \mathcal{R}^u(s):=\sum_{s' \in S_{\mathcal{M},k}}\sum_{a\in\mathcal{A}}\overline{P}^{u} (s'| s,a)\mathcal{R}(s,a)
\end{align}
where $\overline{P}^{u}$$:$$S_{\mathcal{M},k}$$\times$$ \mathcal{A}$$\rightarrow$$\Delta(S_{\mathcal{M},k})$ is defined by replacing parameters $\gamma_a^{q,z}$ and $\delta^{q,z,a}_{q'}$ in \eqref{parametrictransition} with their corresponding instantiations $u(\gamma_a^{q,z})$ and $u(\delta^{q,z,a}_{q'})$. Then, the nonlinear optimization problem to compute the maximum entropy of $\mathcal{D}_{\mathcal{M},k}$ over $\overline{\mathcal{F}}_k(\mathcal{M})$ subject to an expected total reward constraint is
\begin{subequations}
    \begin{align}
       & \underset{\nu,u, \eta}{\text{maximize}} \qquad  \nu(s_{I,\mathcal{M},k})&& \label{objectiveFunction} \\
    &    \text{subject to:}  &&\nonumber \\
       &  \nu(s) \leq L^{u}(s) + \sum_{s' \in S_{\mathcal{M},k}}P_{\mathcal{M},k}^{u} (s'| s)\nu(s') \ \forall \, s \in S_{\mathcal{M},k}\backslash T  \label{fixedPointMapping} \raisetag{10pt}  &&\\
      & \nu(s)=0 \qquad \qquad \qquad \qquad \qquad \qquad \quad\ \  \forall \, s \in T && \label{fixedPointMapping2}\\
      &  \eta(s) \leq \mathcal{R}^u(s)+\sum_{s' \in S_{\mathcal{M},k}}P_{\mathcal{M},k}^{u} (s'| s)\eta(s') \ \forall s \in S_{\mathcal{M},k}  \label{reachabilityProb}  \raisetag{23pt} &&\\
        &  \eta(s_{I,\mathcal{M},k}) \geq \Gamma  \label{initialStateProb}&&\\
       & \sum_{s' \in S_{\mathcal{M},k}}P_{\mathcal{M},k}^{u} (s'| s)=1 \quad \qquad \forall s\in S_{\mathcal{M},k}&& \label{welldefined1}\\
               & P_{\mathcal{M},k}^{u} (s'| s)\geq 0 \quad \quad  \ \qquad  \forall s,s'\in S_{\mathcal{M},k}.&&\label{welldefined2}
    \end{align}
\end{subequations}

As previously explained, the constraints in \eqref{fixedPointMapping}-\eqref{fixedPointMapping2} describe a subspace in $\mathbb{R}^{\lvert S_{\mathcal{M},k}\rvert}$ such that the maximum point $\nu(s_{I,\mathcal{M},k})$ of the subspace corresponds to the value of the maximum entropy $H(\mathcal{D}_{\mathcal{M},k})$ of the pMC $\mathcal{D}_{\mathcal{M},k}$. The constraints \eqref{reachabilityProb}-\eqref{initialStateProb} ensure that the instantiation $u$ satisfies the expected reward constraint given in \eqref{constraint3}. Finally, the constraints \eqref{welldefined1}-\eqref{welldefined2} guarantee that the optimization is performed only over well-defined instantiations $u$.

Note that in the above optimization problem, $P_{\mathcal{M},k}^{u} (s'| s)$, $\eta(s')$ and $\nu(s)$ are functions of decision variables. Therefore, the constraints \eqref{fixedPointMapping} and \eqref{reachabilityProb} contain bilinear terms. Additionally, $L^{u}(s)$ is concave in $P_{\mathcal{M},k}^{u} (s'| s)$, and $\mathcal{R}^u(s)$ is affine in $u(\gamma^{q,z}_a)$ since we consider $k$-FSCs with fixed memory transitions \eqref{LastLoopDetFSC}. 

To solve the optimization problem \eqref{objectiveFunction}-\eqref{welldefined2}, we use a variation of convex-concave-procedure (CCP) \cite{Yuille}, called \textit{penalty} CCP \cite{Lipp}. In particular, we utilize the parameter synthesis method explained in \cite{Cubuktepe:10.1007/978-3-030-01090-4_10}. Here, we briefly explain the solution approach and refer the reader to \cite{Cubuktepe:10.1007/978-3-030-01090-4_10} for details. 

We first represent each bilinear term $f(x)$, e.g., $P_{\mathcal{M},k}^{u} (s'| s)\nu(s')$, as a difference-of-convex function $f(x)$$=$$f_1(x)$$-$$f_2(x)$ and linearize the concave part $f_2(x)$ around an initial point. Doing so yields a convex optimization problem. We then introduce nonnegative penalty variables $\psi_i$ to the constraints \eqref{fixedPointMapping} and \eqref{reachabilityProb}, and replace the objective function with $\nu(s_{I,\mathcal{M},k})$$-$$\tau\sum_i\psi_i$ where $\tau$ is a constant regularization parameter. We solve the resulting convex problem and update the initial point with the optimal solution to the convex problem. By iteratively performing the same steps, we obtain, if the procedure converges, a local optimal solution to our original problem  \eqref{objectiveFunction}-\eqref{welldefined2}.

\section{Numerical Examples}
We now provide two numerical examples to demonstrate the relation between the maximum entropy of a POMDP, the threshold $\Gamma$ on the expected total reward, and the number of memory states in the FSCs. We use the MOSEK \cite{mosek} solver with the CVX \cite{cvx} interface to solve the convex optimization problems obtained from the convex-concave procedure. To improve the approximation of exponential cone constraints, we use the CVXQUAD \cite{cvxquad} package. 

\subsection{Relation Between the Maximum Entropy and the Expected Reward Threshold}
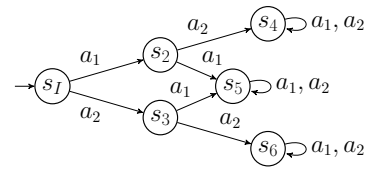
\begin{figure}[b!]
    \centering
    \scalebox{0.55}{
        \begin{tikzpicture}[->, >=stealth', auto, semithick, node distance=2cm]
            \tikzstyle{every state}=[fill=white,draw=black,thick,text=black,scale=0.7]

            \node[state,initial,initial text=] (s_0) {\Huge{$s_I$}};
            \node[state] (s_1) [above right = 1.5mm and 20mm of s_0] {\Huge{$s_2$}};
            \node[state] (s_2) [below right = 1.5mm and 20mm of s_0] {\Huge{$s_3$}};
            \node[state] (s_3) [above right = 1.5mm and 20mm of s_1] {\Huge{$s_4$}};
            \node[state] (s_4) [right = 35mm of s_0] {\Huge{$s_5$}};
            \node[state] (s_5) [below right = 1.5mm and 20mm of s_2] {\Huge{$s_6$}};
        
            \path (s_0) edge node[above left]{\LARGE{$a_1$}} (s_1)
                  (s_0) edge node[below left]{\LARGE{$a_2$}} (s_2)
                  (s_1) edge node[above left]{\LARGE{$a_2$}} (s_3)
                  (s_1) edge node[above right]{\LARGE{$a_1$}} (s_4)
                  (s_2) edge node[above left]{\LARGE{$a_1$}} (s_4)
                  (s_2) edge node[above right]{\LARGE{$a_2$}} (s_5)
                  (s_3) edge[loop right] node[right]{\LARGE{$a_1,a_2$}} (s_3)
                  (s_4) edge[loop right] node[right]{\LARGE{$a_1,a_2$}} (s_4)
                  (s_5) edge[loop right] node[right]{\LARGE{$a_1,a_2$}} (s_5);
                  
        \end{tikzpicture}}
    \caption{POMDP illustrating the relation between the maximum entropy and the expected total reward $\Gamma$.}
    \label{fig:simpleExample}
\end{figure}

In the first example, we consider a POMDP with 6 states shown in Fig. \ref{fig:simpleExample}. There is only one observation $\mathcal{Z}=\{z_1\}$ and therefore the observation function is $\mathcal{O}_{s,z_1}$$=$$1$ for all states \textit{s}. We use a deterministic 2-FSC whose memory transition function $\delta$ is given in \eqref{LastLoopDetFSC}. Because there is only one observation, the synthesized controller is an open-loop controller. We suppose that the agent aims to reach state $s_{4}$ and encode this objective by defining a reward function $\mathcal{R}$ such that $\mathcal{R}(s_2,a_1)$$=$$\mathcal{R}(s_3,a_1)$$=$$1$ and $\mathcal{R}(s,a)$$=$$0$ otherwise.

\begin{figure}[t!]
    \centering
    \scalebox{0.6}{
    \begin{tikzpicture}
        \begin{axis}[
            xlabel={\Large{Expected Total Reward Threshold ($\Gamma$)}},
            ylabel={\Large{Entropy of Stochastic Process [bits]}},
            xmin=0.5,xmax=1,
            ymin=0.70,ymax=2.05,
            xtick={0.5,0.6,0.7,0.8,0.9,1},
            ytick={0.75,1.00,1.25,1.50,1.75,2.00},
            legend pos=south west,
            ymajorgrids=true,
            grid style=dashed]
            \addplot[color=blue,mark=square]
                coordinates{(0.5,2.00)(0.6,1.97)(0.7,1.88)(0.8,1.72)(0.9,1.47)(1.0,1.00)};
            \addplot[color=green,mark=triangle]
                coordinates{(0.5,1.9750)(0.6,1.9235)(0.7,1.8209)(0.8,1.6540)(0.9,1.416)(1.0,0.8367)};
            \addplot[color=red,mark=*,dashed]
                coordinates{(0.5,2.00)(0.6,1.97)(0.7,1.88)(0.8,1.72)(0.9,1.47)(1.0,1.00)};
            \legend{\large{Proposed method}, \large{Feasibility approach}, \large{MDP upper bound}}
        \end{axis}
    \end{tikzpicture}}
    \caption{The trade-off between the maximum entropy and the expected total rewards.}
    \label{fig:EntMaxVsFeas}
\end{figure}
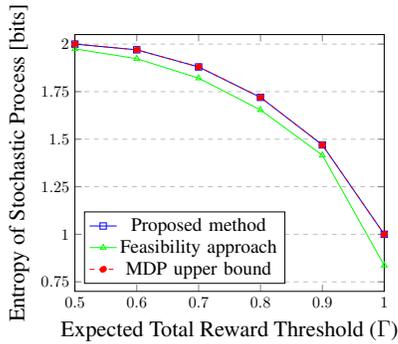

We investigate the effect of the threshold $\Gamma$ in \eqref{constraint3} on the maximum entropy by synthesizing controllers for values between $\Gamma$$=$$0.5$ and $\Gamma$$=$$1$. For each value of $\Gamma$, we solve the optimization problem given in Section \ref{FSCsynthesis} for 10 times by randomly initializing the convex-concave procedure. For each $\Gamma$, we pick the best result of 10 trials, and plot the maximum entropy of the stochastic process induced by the synthesized controllers in Fig. \ref{fig:EntMaxVsFeas}. For comparison, we synthesize controllers by solving a feasibility problem given in \cite{Cubuktepe:10.1007/978-3-030-01090-4_10}. We obtain the feasibility problem from \eqref{objectiveFunction}-\eqref{welldefined2} by removing the entropy constraint (\ref{fixedPointMapping}) and replacing the objective function (\ref{objectiveFunction}) with a constant value.

In this example, the proposed approach yields the globally optimal controller by attaining a tight bound on $\Gamma$. The global optimality of the controller is evident in Figure \ref{fig:EntMaxVsFeas}, as the entropy of the proposed approach exactly matches that of the underlying MDP for each value of $\Gamma$. Because the feasibility program only seeks to find a feasible instantiation of the parameters that satisfy the expected total reward constraint (\ref{reachabilityProb}), the entropy of the stochastic processes it yields is less than the maximum attainable entropy.

\subsection{Relation Between the Maximum Entropy and the Number of Memory States}

\begin{figure}[t!]
    \centering\scalebox{0.65}{
    \begin{tikzpicture}[->, >=stealth', auto, semithick, node distance=2cm]
        \tikzstyle{every state}=[fill=white,draw=black,thick,text=black,scale=0.7]

        \node[state,initial,initial text=] (s_2) {\Huge{$s_I$}};
        \node[state] (s_1) [above = 5mm of s_2] {\Huge{$s_1$}};
        \node[state] (s_3) [below = 5mm of s_2] {\Huge{$s_3$}};
        
        \node[state] (s_4) [right = 8mm of s_1] {\Huge{$s_4$}};
        \node[state] (s_7) [right = 8mm of s_4] {\Huge{$s_7$}};
        \node[state] (s_10) [right = 8mm of s_7] {\Huge{$s_{10}$}};
        \node[state] (s_13) [right = 8mm of s_10] {\Huge{$s_{13}$}};
        
        \node[state] (s_5) [right = 8mm of s_2] {\Huge{$s_5$}};
        \node[state] (s_8) [right = 8mm of s_5] {\Huge{$s_8$}};
        \node[state] (s_11) [right = 8mm of s_8] {\Huge{$s_{11}$}};
        \node[state] (s_14) [right = 8mm of s_11] {\Huge{$s_{14}$}};
        
        \node[state] (s_6) [right = 8mm of s_3] {\Huge{$s_6$}};
        \node[state] (s_9) [right = 8mm of s_6] {\Huge{$s_9$}};
        \node[state] (s_12) [right = 8mm of s_9] {\Huge{$s_{12}$}};
        \node[state] (s_15) [right = 8mm of s_12] {\Huge{$s_{15}$}};
        
        \node (dots_1) [right = 1mm of s_2] {$\cdots$};
        \node (dots_2) [right = 27mm of dots_1] {$\cdots$};
        \node (dots_3) [right = 10mm of dots_2] {$\cdots$};
        
        \node (dots_4) [right = 1mm of s_4] {$\cdots$};
        \node (dots_5) [right = 10mm of dots_4] {$\cdots$};
        \node (dots_6) [right = 10mm of dots_5] {$\cdots$};
        
        \node (dots_7) [right = 1mm of s_6] {$\cdots$};
        \node (dots_8) [right = 10mm of dots_7] {$\cdots$};
        \node (dots_9) [right = 10mm of dots_8] {$\cdots$};
        
        \path (s_13) edge[loop right] node[right]{} (s_13)
            (s_14) edge[loop right] node[right]{} (s_14)
            (s_15) edge[loop right] node[right]{} (s_15);
        
        \path (s_1) edge node[above]{\large{$a_1$}} (s_4)
            (s_1) edge node[above]{\large{$a_2$}} (s_5)
            (s_1) edge node[left]{\large{$a_3$}} (s_2);
            
        \path (s_5) edge node[above]{\large{$a_1$}} (s_7)
            (s_5) edge node[above]{\large{$a_2$}} (s_8)
            (s_5) edge node[above]{\large{$a_3$}} (s_9);
            
        \path (s_3) edge node[above]{\large{$a_1$}} (s_6)
            (s_3) edge node[above]{\large{$a_2$}} (s_5)
            (s_3) edge node[left]{\large{$a_3$}} (s_2);
             
        \end{tikzpicture}}
    \caption{POMDP illustrating the relation between the maximum entropy and the number of memory states in FSCs.}
    \label{fig:POMDPex2}
\end{figure}
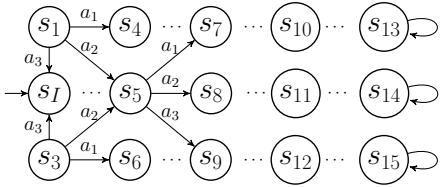

In this example, we consider a POMDP with 15 states shown in Fig. \ref{fig:POMDPex2}. As in the previous example, there is only a single observation $\mathcal{Z}$$=$$\{z_{1}\}$ yielding $\mathcal{O}_{s,z_{1}}$$=$$1$ for all states $s$. We suppose that the agent aims to reach $s_{14}$ with probability 1. To encode this objective, we set $\Gamma$$=$$1$ with $\mathcal{R}(s_{10},a_{2})$$=$$1$, $\mathcal{R}(s_{11},a_{2})$$=$$1$, $\mathcal{R}(s_{12},a_{2})$$=$$1$, and $\mathcal{R}(s,a)$$=$$0$ otherwise.

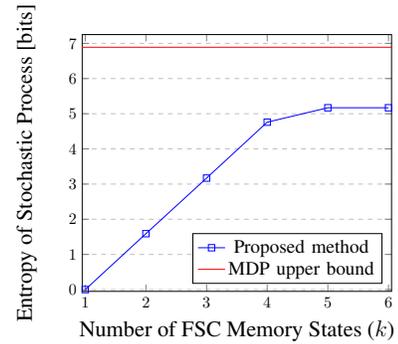
\begin{figure}[t]
    \centering
    \scalebox{0.6}{
    \begin{tikzpicture}
        \begin{axis}[
            xlabel={\Large{Number of FSC Memory States ($k$)}},
            ylabel={\Large{Entropy of Stochastic Process [bits]}},
            xmin=0.95,xmax=6.05,
            ymin=-0.05,ymax=7.25,
            xtick={0,1,2,3,4,5,6},
            ytick={0,1.0,2.0,3.0,4.0,5.0,6.0,7.0},
            ymajorgrids=true,
            grid style=dashed,
            legend pos=south east]
            \addplot[color=blue,mark=square]
                coordinates{(1,0)(2,1.59)(3,3.17)(4,4.76)(5,5.17)(6,5.17)};
            \addplot[color=red,mark=none]
                coordinates{(0,6.89)(2,6.89)(3,6.89)(4,6.89)(5,6.89)(7,6.89)};
            \legend{\large{Proposed method}, \large{MDP upper bound}};
        \end{axis}
    \end{tikzpicture}}
    \caption{Comparison between the maximum entropy of the induced stochastic process for varying values of $k$.}
    \label{fig:EntropyVsk}
\end{figure}

\begin{figure}[t!]
    \centering\scalebox{0.6}{
        \begin{tikzpicture}[->, >=stealth', auto, semithick, node distance=2cm]
        \tikzstyle{every state}=[fill=white,draw=black,thick,text=black,scale=0.7]

        \node[state,initial,initial text=] (s_2) {\Huge{$s_I$}};
        \node[state] (s_1) [above = 5mm of s_2] {\Huge{$s_1$}};
        \node[state] (s_3) [below = 5mm of s_2] {\Huge{$s_3$}};
        
        \node[state] (s_4) [right = 8mm of s_1] {\Huge{$s_4$}};
        \node[state] (s_7) [right = 8mm of s_4] {\Huge{$s_7$}};
        \node[state] (s_10) [right = 8mm of s_7] {\Huge{$s_{10}$}};
        \node[state] (s_13) [right = 8mm of s_10] {\Huge{$s_{13}$}};
        
        \node[state] (s_5) [right = 8mm of s_2] {\Huge{$s_5$}};
        \node[state] (s_8) [right = 8mm of s_5] {\Huge{$s_8$}};
        \node[state] (s_11) [right = 8mm of s_8] {\Huge{$s_{11}$}};
        \node[state] (s_14) [right = 8mm of s_11] {\Huge{$s_{14}$}};
        
        \node[state] (s_6) [right = 8mm of s_3] {\Huge{$s_6$}};
        \node[state] (s_9) [right = 8mm of s_6] {\Huge{$s_9$}};
        \node[state] (s_12) [right = 8mm of s_9] {\Huge{$s_{12}$}};
        \node[state] (s_15) [right = 8mm of s_12] {\Huge{$s_{15}$}};
        
        \draw[line width=0.15mm, red] (s_2) -- (s_4);
        \draw[line width=0.15mm, red] (s_2) -- (s_5);
        \draw[line width=0.15mm, red] (s_2) -- (s_6);
        
        \draw[line width=0.15mm, red] (s_4) -- (s_5);
        \draw[line width=0.15mm, red] (s_4) -- (s_7);
        \draw[line width=0.15mm, red] (s_4) -- (s_8);
        
        \draw[line width=0.15mm, red] (s_5) -- (s_7);
        \draw[line width=0.15mm, red] (s_5) -- (s_8);
        \draw[line width=0.15mm, red] (s_5) -- (s_9);
        
        \draw[line width=0.15mm, red] (s_6) -- (s_5);
        \draw[line width=0.15mm, red] (s_6) -- (s_8);
        \draw[line width=0.15mm, red] (s_6) -- (s_9);
        
        \draw[line width=0.15mm, red] (s_7) -- (s_8);
        \draw[line width=0.15mm, red] (s_7) -- (s_10);
        \draw[line width=0.15mm, red] (s_7) -- (s_11);
        
        \draw[line width=0.15mm, red] (s_8) -- (s_10);
        \draw[line width=0.15mm, red] (s_8) -- (s_11);
        \draw[line width=0.15mm, red] (s_8) -- (s_12);
        
        \draw[line width=0.15mm, red] (s_9) -- (s_8);
        \draw[line width=0.15mm, red] (s_9) -- (s_11);
        \draw[line width=0.15mm, red] (s_9) -- (s_12);
        
        \draw[line width=0.30mm, red] (s_10) -- (s_14);
        \draw[line width=0.30mm, red] (s_10) -- (s_11);
        
        \draw[line width=0.45mm, red] (s_11) -- (s_14);
        
        \draw[line width=0.45mm, red] (s_12) -- (s_14);
        \draw[line width=0.30mm, red] (s_12) -- (s_11);
        
        \path[line width=0.45mm, red] (s_14) edge [out=45,in=315,looseness=6] node[above] {} (s_14);
        \end{tikzpicture}}
    \caption{Trajectories of synthesized entropy maximizing controller. Edge thicknesses indicate the transition probabilities.}
    \label{fig:POMDPoptController}
\end{figure}
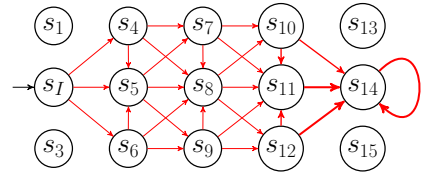

We study the relation between the number of memory states and the maximum entropy of the induced pMC by synthesizing controllers for $k$$=$$1,\ldots,6$ memory states. As in the previous example, we run the optimization problem given in Section \ref{FSCsynthesis} for 10 times while randomly initializing the convex-concave procedure. In Fig. \ref{fig:EntropyVsk}, we plot the maximum entropy of the stochastic process induced by the controller for each value of $k$. Furthermore, Fig. \ref{fig:POMDPoptController} shows the entropy-maximizing controller for the POMDP, where edge weights correspond to the probability of action selection. 

From Fig. \ref{fig:EntropyVsk}, we see that the 1-FSC achieves a maximum entropy of 0. A 1-FSC selecting any action besides $a_{2}$ cannot reach state $s_{14}$ while collecting an expected total reward of 1. An additional memory state allows the agent to randomize its action selection for one more time step. After 5 memory states, however, additional memory states do not affect the maximum entropy of the induced stochastic process. 

Any controller with at least 5 memory states achieves an optimal action distribution shown in Fig. $\ref{fig:POMDPoptController}$. Unlike the previous example, a gap between the maximum entropy of the MDP and that of the induced pMC remains. The maximum entropy of the POMDP must lie within this gap. This example demonstrates the monotonicity of the maximum entropy with the number of states in the controller.

\section{Conclusions}
In this paper, we consider an entropy maximization problem in POMDPs subject to an expected reward constraint. We define the entropy in POMDPs and show that the maximum entropy of a POMDP is upper bounded by that of the underlying MDP. We then consider the entropy maximization problem over deterministic FSCs. We show that this problem can be translated to the so-called parameter synthesis problem in a pMC obtained by the product of the POMDP and the FSC. We propose to use penalty CCP to solve such a nonlinear optimization problem. Two examples are presented to show the validity of our proposed approach. 

\bibliographystyle{IEEEtran}
\bibliography{ref1}

\section{Appendix}

\textbf{Proof of Lemma \ref{writeValueFunction}.} We prove the claim by strong induction on \textit{t}. For the base case, we have

\begin{subequations}
    \allowdisplaybreaks
    \begin{align}
       & \mathcal{V}_{T-1,T}^{\pi} (s^{T-1})= H^{\pi}(X_{T}|X_{T-1},X^{T-1}=s^{T-1}) \nonumber \\
        & \quad + H^{\pi}(X_{T+1}|X_{T},X^{T-1}=s^{T-1}) \label{FirstInd-a} \\
        & = H^{\pi}(X_{T}|X^{T-1}=s^{T-1}) \nonumber \\
        & \quad  + H^{\pi}(X_{T+1}|X_{T},\mathcal{H}^{T-1}=h^{T-1}) \label{FirstInd-b} \\
        & = H^{\pi}(X_{T}|X^{T-1}=s^{T-1}) \nonumber \\
        &\quad + \sum_{s^{T} \in \mathcal{S}\mathcal{H}^{T}} Pr^{\pi}(s^{T}|s^{T-1}) H^{\pi}(X_{T+1}|X^{T}=s^{T}) \label{FirstInd-d} \\
        & =  H^{\pi}(X_{T}|X^{T-1}=s^{T-1}) \nonumber \\
        &\quad + \sum_{s^{T} \in \mathcal{S}\mathcal{H}^{T}} Pr^{\pi}(s^{T}|s^{T-1}) \mathcal{V}_{T,T}^{\pi}(s^{T}) \label{FirstInd-e}.
    \end{align}
\end{subequations}

\noindent
where (\ref{FirstInd-b}) follows from (\ref{FirstInd-a}) by the fact that $X_{T-1}$ is a component of $S^{T-1}$. By the total law of probability and the definition of the state history, we obtain (\ref{FirstInd-d}) from (\ref{FirstInd-b}). Lastly, (\ref{FirstInd-e}) holds by the definition of the value function defined in (\ref{eq:EntropyVfRecursive}). We now assume that the equality in (\ref{eq:EntropyVfRecursive}) holds for time steps $T-2,T-3,\ldots,t+1$, and show that the equality holds for $t$.
\begin{subequations}
\allowdisplaybreaks
    \begin{align}
        \mathcal{V}_{t,T}^{\pi}(s^{t}) &  = \sum_{k=t}^{T} H^{\pi}(X_{k+1}|X_{t}^{k},X^{t}=s^{t}) \label{SecondInd-a} \\
        & = H^{\pi}(X_{t+1}|X_{t},X^{t}=s^{t}) \nonumber \\
        & \qquad + \sum_{k=t+1}^{T} H^{\pi}(X_{k+1}|X_{t+1}^{k},X^{t}=s^{t}) \label{SecondInd-b} \\
        & = H^{\pi}(X_{t+1}|X^{t}=s^{t}) + \sum_{s^{t+1} \in \mathcal{S}\mathcal{H}^{t}} \sum_{k=t+1}^{T}... \nonumber \\
        & \qquad H^{\pi} Pr^{\pi}(s^{t+1}|s^{t}) (X_{k+1}|X_{t}^{k},X^{t+1}=s^{t+1}) \label{SecondInd-d} \\
        & = H^{\pi}(X_{t+1}|X^{t}=s^{t}) \nonumber \\
        & \qquad + \sum_{s^{t+1} \in \mathcal{S} \mathcal{H}^{t+1}} Pr^{\pi}(s^{t+1}|s^{t}) \mathcal{V}_{t+1,T}^{\pi}(s^{t+1}). \label{SecondInd-e}
    \end{align}
\end{subequations}
As in the base case, (\ref{SecondInd-b}) follows from (\ref{SecondInd-a}) by the fact that $X_{t}$ is a component of $s^{t}$. We then obtain (\ref{SecondInd-d}) from (\ref{SecondInd-b}) by the total law of probability and the definition of the state history $h^{t}$. Lastly, (\ref{SecondInd-e}) holds by the definition of the value function defined in (\ref{eq:EntropyVfRecursive}). The equality holds for a general $t$, completing the induction. We may thus write the total expected entropy in this recursive form.$\Box$

\textbf{Proof of Theorem \ref{thm:POMDPbounded}.} We prove the claim by strong induction on $t$. Denote the value function for $\pi \in \Pi(\mathcal{M})$ as $\mathcal{V}_{t,T}^{\pi}(s^{t})$ and the value function for $\pi'$$\in$$\Pi(\mathcal{M}_{fo})$ constructed according to (\ref{MDP_controller}) as $\mathcal{V}_{t,T}^{\pi'}(s^{t})$, respectively. Starting with the base case $t=T$,  we have
\begin{subequations}
    \begin{align}
        \mathcal{V}_{T,T}^{\pi}(s^{T}) & = H^{\pi}(X_{T+1}|X^{T}=s^{T}) \label{eq:InductLastStep1-a} \\
        & = H^{\pi'}(X_{T+1}|X^{T}=s^{T}) \label{eq:InductLastStep1-b} \\
        \sup_{\pi \in \Pi(\mathcal{M})} \mathcal{V}_{T,T}^{\pi}(s^{T}) & \leq \sup_{\pi' \in \Pi(\mathcal{M}_{fo})} H^{\pi'}(X_{T+1}|X^{T}=s^{T}) \label{eq:InductLastStep1-c} \\
        & = \mathcal{V}_{T,T}^{\pi'}(s^{T}). \label{eq:InductLastStep1-d}
    \end{align}
\end{subequations}
The equality in (\ref{eq:InductLastStep1-b}) follows from the fact that we can construct an equivalent history-dependent controller on the underlying MDP that achieves the same transition probabilities for any observation-based controller. We then obtain \eqref{eq:InductLastStep1-c} by the fact that $\Pi(\mathcal{M})$$\subset$$\Pi(\mathcal{M}_{fo})$. By the definition of the value function in (\ref{eq:EntropyVfRecursive}), we then obtain (\ref{eq:InductLastStep1-c}). 

Now assume that the inequality holds for time steps $T$$-$$1, ..., t+1$. We show that it also holds for $t$ as follows. Note first that
\begin{subequations}
    \begin{align}
        &\mathcal{V}_{t,T}^{\pi}(s^{t}) = H^{\pi}(X_{t+1}|X^{t}=s^{t}) \nonumber \\
        & \qquad + \sum_{\substack{s^{t+1} \in \mathcal{S}\mathcal{H}^{t+1}}} Pr^{\pi}(s^{t+1}|s^{t}) \mathcal{V}_{t+1,T}^{\pi}(s^{t+1}) \label{eq:InductMidStep1-a} \\
        &\qquad \leq H^{\pi}(X_{t+1}|X^{t}=s^{t}) \nonumber \\
        & \qquad + \sum_{\substack{s^{t+1} \in \mathcal{S}\mathcal{H}^{t+1}}} Pr^{\pi}(s^{t+1}|s^{t}) \mathcal{V}_{t+1,T}^{\pi'}(s^{t+1}) \label{eq:InductMidStep1-b} \\
        &\qquad = H^{\pi'}(X_{t+1}|X^{t}=s^{t}) \nonumber \\
        & \qquad + \sum_{\substack{s^{t+1} \in \\ \mathcal{S}\mathcal{H}^{t+1}}} \mathcal{V}_{t+1,T}^{\pi'} (s^{t+1}) Pr^{\pi'}(s^{t+1}|s^{t+1}). \label{eq:InductMidStep1-c} 
    \end{align}
\end{subequations}
By Lemma \ref{writeValueFunction}, we can write the value function recursively in (\ref{eq:InductMidStep1-a}). The equality in (\ref{eq:InductMidStep1-b}) then follows by the induction hypothesis. By (\ref{MDP_controller}), we can construct an equivalent controller on the underlying MDP that has the same transition probabilities. Doing so yields (\ref{eq:InductMidStep1-c}). Then, we have
\begin{subequations}
\begin{align}
            \sup_{\pi \in \Pi(\mathcal{M})}\mathcal{V}_{t,T}^{\pi}(s^{t}) & \leq \sup_{\pi' \in \Pi(\mathcal{M}_{fo})} H^{\pi'}(X_{t+1}|X^{t}=s^{t}) \nonumber \\
        & \, + \sum_{\substack{s^{t+1} \in \\ \mathcal{S}\mathcal{H}^{t+1})}} Pr^{\pi'}(s^{t+1}|s^{t}) \mathcal{V}_{t+1,T}^{\pi'}(s^{t+1}) \label{eq:InductMidStep1-d} \\
        & = \mathcal{V}_{t,T}^{\pi'}(s^{t}). \label{eq:InductMidStep1-e}
\end{align}
\end{subequations}
where inequality in (\ref{eq:InductMidStep1-d}) is due to the fact that $\Pi(\mathcal{M})$$\subset$$\Pi(\mathcal{M}_{fo})$ and (\ref{eq:InductMidStep1-e}) follows by the definition of the value function in (\ref{eq:EntropyVfRecursive}). Thus the induction holds for $t$. Since the claim holds for all $t$, we have $\mathcal{V}_{1,T}^{\pi}(s_{I})$$\leq$$\mathcal{V}_{1,T}^{\pi'}(s_{I})$. By (\ref{eq:defEntropyMDP}), this implies that $H^{\pi}(X^{T})$$\leq$$H^{\pi'}(X^{T})$ for all $T$. Taking the limit as $T$$\rightarrow$$\infty$ on both sides of the inequality completes the proof. $\Box$

\textbf{Proof of Theorem \ref{thm:pMCbounded}.} By definition of the $H(\mathcal{D}_{\mathcal{M},k})$, each possible instantiation $u_{\mathcal{C}}$ can only correspond to a deterministic FSC $\mathcal{C}$, i.e., all corresponding FSCs satisfy $\lvert Succ(q)\rvert$$=$$1$. Then, it can be shown by construction that there is a one-to-one correspondence between the state histories of the instantiated pMC $\mathcal{D}_{\mathcal{M},k}[u_{\mathcal{C}}]$ and its corresponding POMDP $\mathcal{M}$ under the FSC $\mathcal{C}$. Additionally, for any given instantiation $u_{\mathcal{C}}$, there exists a deterministic FSC $\mathcal{C}$$\in$$\mathcal{F}_{k}^{det}(\mathcal{M})$ which induces the same state history transition function $Pr^{u_{\mathcal{C}}}(s^{t+1}| s^t)$ with $\mathcal{D}_{\mathcal{M},k}[u_{\mathcal{C}}]$. Therefore, using the result of Lemma $\ref{writeValueFunction}$, we can show that for any instantiated pMC, there exists an FSC that will induce from the POMDP a stochastic process with the same entropy. Because $\mathcal{F}_{k}^{det}(\mathcal{M})$$\subset$$\Pi(\mathcal{M})$, it then follows that $H(\mathcal{D}_{\mathcal{M},k})$$\leq$$H(\mathcal{M})$.$\Box$

\textbf{Proof of Lemma \ref{lemma:pMCmonotonic}.} We prove the claim by induction on the number of memory states $k$. We start with the base case $n$$=$$1$. Consider an instantiated pMC $\mathcal{D}_{\mathcal{M},1}[u_{\mathcal{C}}]$ for which there exists a corresponding deterministic 1-FSC $\mathcal{C}$$\in$$\bar{\mathcal{F}}_{1}(\mathcal{M})$ whose decision function $\gamma$ satisfies $\gamma(a|q_{1},z)$$=$$u_{\mathcal{C}}(\gamma_{a}^{q_{1},z})$. Now, construct a deterministic 2-FSC $\mathcal{C}'$ whose decision function $\gamma'$ satisfies $\gamma'(a|q_{1},z)$$=$$\gamma'(a|q_{2},z)$$=$$u_{\mathcal{C}}(\gamma_{a}^{q_{1},z})$. Then, since the memory transitions of both FSCs satisfy \eqref{LastLoopDetFSC}, there is a one to one correspondence between the state histories of $\mathcal{D}_{\mathcal{M},1}[u_{\mathcal{C}}]$ and $\mathcal{D}_{\mathcal{M},2}[u_{\mathcal{C}'}]$. Using Lemma \ref{writeValueFunction}, it can be shown that $H(\mathcal{D}_{\mathcal{M},1}[u_{\mathcal{C}}])$$=$$H(\mathcal{D}_{\mathcal{M},2}[u_{\mathcal{C}'}])$. Since we choose $\mathcal{C}$ arbitrarily, the maximum entropy of $\mathcal{D}_{\mathcal{M},2}$ cannot be lower than that of $\mathcal{D}_{\mathcal{M},1}$, i.e.,
\begin{equation}
    \sup_{\mathcal{C} \in \overline{\mathcal{F}}_{1}(\mathcal{M})}H(\mathcal{D}_{\mathcal{M},1}[u_{\mathcal{C}}])\leq \sup_{\mathcal{C} \in \overline{\mathcal{F}}_{2}(\mathcal{M})}H(\mathcal{D}_{\mathcal{M},2}[u_{\mathcal{C}}]).
\end{equation}

We assume that the claim holds for $n$$=$$1,2,...,k-1$, and show that it also holds for $k$$=$$n$. Consider an instantiated pMC $\mathcal{D}_{\mathcal{M},k-1}[u_{\mathcal{C}}]$ for which there exists a corresponding deterministic $(k-1)$-FSC $\mathcal{C}$$\in$$\bar{\mathcal{F}}_{k-1}(\mathcal{M})$ whose decision function $\gamma$ satisfies $\gamma(a|q_{i},z)$$=$$u_{\mathcal{C}}(\gamma_{a}^{q_{i},z})$ for $i$$=$$1,\ldots,k-1$. Then, we can construct an $k$-FSC $\mathcal{C}'$ whose decision function $\gamma'$ satisfies $\gamma'(a|q_{i},z):=\gamma(a|q_{i},z)$ for $i$$=$$1,\ldots,k-1$, and $\gamma'(a|q_{k},z):=\gamma(a|q_{k-1},z)$. Since the memory transitions of both FSCs satisfy \eqref{LastLoopDetFSC}, there is a one to one correspondence between the state histories of $\mathcal{D}_{\mathcal{M},k-1}[u_{\mathcal{C}}]$ and $\mathcal{D}_{\mathcal{M},k}[u_{\mathcal{C}'}]$. Using Lemma \ref{writeValueFunction}, it can be shown that $H(\mathcal{D}_{\mathcal{M},k-1}[u_{\mathcal{C}}])$$=$$H(\mathcal{D}_{\mathcal{M},k}[u_{\mathcal{C}'}])$. Then, since $\mathcal{C}$ is chosen arbitrarily, using the induction hypothesis, we obtain
\begin{equation}
    \sup_{\mathcal{C} \in \overline{\mathcal{F}}_{j}(\mathcal{M})}H(\mathcal{D}_{\mathcal{M},j}[u_{\mathcal{C}}])\leq \sup_{\mathcal{C} \in \overline{\mathcal{F}}_{k}(\mathcal{M})}H(\mathcal{D}_{\mathcal{M},k}[u_{\mathcal{C}}])
\end{equation}  
for all $j$$\leq$$k$. This completes the proof. $\Box$

\end{document}